\documentclass[14pt]{article}
\usepackage{mathrsfs}
\usepackage{amsthm}
\usepackage{amssymb}
\usepackage{amsmath}
\usepackage{graphicx}
\usepackage{color}
\usepackage{amsfonts}
\usepackage{float}
\usepackage{cite}
\usepackage[text={140mm,210mm},left=45mm,vmarginratio=1:1]{geometry}
\newtheorem{theorem}{Theorem}[section]

\newtheorem{lemma}[theorem]{Lemma}
\newtheorem{proposition}[theorem]{Proposition}

\newtheorem{corollary}[theorem]{Corollary}

\numberwithin{equation}{section}
\normalsize

\begin{document}
\title{\textbf{Law of large numbers for the SIR model with random vertex weights on Erd\H{o}s-R\'{e}nyi graph}}

\author{Xiaofeng Xue \thanks{\textbf{E-mail}: xfxue@bjtu.edu.cn \textbf{Address}: School of Science, Beijing Jiaotong University, Beijing 100044, China.}\\ Beijing Jiaotong University}

\date{}
\maketitle

\noindent {\bf Abstract:}
In this paper we are concerned with the SIR model with random vertex weights on Erd\H{o}s-R\'{e}nyi graph $G(n,p)$. The Erd\H{o}s-R\'{e}nyi graph $G(n,p)$ is generated from the complete graph $C_n$ with $n$ vertices through independently deleting each edge with probability $(1-p)$. We assign i. i. d. copies of a positive r. v. $\rho$ on each vertex as the vertex weights. For the SIR model, each vertex is in one of the three states `susceptible', `infective' and `removed'. An infective vertex infects a given susceptible neighbor at rate proportional to the production of the weights of these two vertices. An infective vertex becomes removed at a constant rate. A removed vertex will never be infected again. We assume that at $t=0$ there is no removed vertex and the number of infective vertices follows a Bernoulli distribution $B(n,\theta)$. Our main result is a law of large numbers of the model. We give two deterministic functions $H_S(\psi_t), H_V(\psi_t)$ for $t\geq 0$ and show that for any $t\geq 0$, $H_S(\psi_t)$ is the limit proportion of susceptible vertices and $H_V(\psi_t)$ is the limit of the mean capability of an infective vertex to infect a given susceptible neighbor at moment $t$ as $n$ grows to infinity. 

\quad


\section{Introduction}\label{section one}
In this paper we are concerned with the SIR (Susceptible-Infective-Removed) model with random vertex weights on Erd\H{o}s-R\'{e}nyi graph $G(n,p)$. First we introduce some notations. For each integer $n\geq 1$, we denote by $A_n$ the set $\{0,1,2,\ldots,n-1\}$. We consider the $n$ elements in $A_n$ as $n$ vertices and assume that any two vertices are connected by an edge. As a result, we obtain a complete graph with $n$ vertices, which we denote as $C_n$. Let $p\in (0,1)$, then we can obtain a random graph $G_n$ through the procedure that each edge on $C_n$ is independently deleted with probability $1-p$, in other words, remained with probability $p$. The graph $G_n$ with vertices set $A_n$ and edges which are remained is called the Erd\H{o}s-R\'{e}nyi graph with parameter $G(n,p)$ (see Chapter 4 of \cite{Hofstad2013}). For any $0\leq i<j\leq n-1$, we denote by $i\sim j$ when the edge connecting $i$ and $j$ on $C_n$ is remained during the procedure to generate $G_n$. That is to say, $i\sim j$ when and only when $i$ is a neighbor of $j$ on the graph $G_n$.

Let $\rho$ be a positive random variable such that $P(\rho>0)>0$ and $P(0\leq\rho\leq M_1)=1$ for some $M_1<+\infty$ while $\{\rho(i)\}_{0\leq i<+\infty}$ are i. i. d. copies of $\rho$ which are independent of $\{G_n\}_{n\geq 1}$, then the SIR model $\{X_t\}_{t\geq 0}$ on $G_n$ with vertex weights $\{\rho(i)\}_{0\leq i\leq n-1}$ is a continuous-time Markov process with state space $\{0,1,-1\}^{A_n}$. That is to say, at each moment $t\geq 0$, there is a spin on each vertex on $G_n$ with value taking from $\{0,1,-1\}$. For each $i\in A_n$ and $t\geq 0$, we denote by $X_t(i)$ the value of the spin on $i$ at moment $t$, then $\{X_t\}_{t\geq 0}$ evolves according to the following rules. If $X_t(i)=-1$, then $i$ is frozen in state $-1$ after the moment $t$. That is to say, $X_s(i)=-1$ for any $s\geq t$. If $X_{t}(i)=0$, then
\[
P\big(X_{t+\Delta t}(i)=1\big|X_s,s\leq t\big)=\Big(\frac{\lambda}{n}\sum\limits_{j=0}^{n-1}\rho(i)\rho(j)\mathbf{1}_{\{j\sim i, X_t(j)=1\}}\Big)\Delta t+o(\Delta t)
\]
and
\[
P\big(X_{t+\Delta t}(i)=0\big|X_s,s\leq t\big)=1-P\big(X_{t+\Delta t}(i)=1\big|X_s,s\leq t\big)+o(\Delta t)
\]
where $\lambda>0$ is a positive parameter called the infection rate and $\mathbf{1}_A$ is the indicator function of the random event $A$. If $X_{t}(i)=1$, then
\[
P\big(X_{t+\Delta t}(i)=-1\big|X_s,s\leq t\big)=\Delta t+o(\Delta t)=1-P\big(X_{t+\Delta t}(i)=1\big|X_s,s\leq t\big)+o(\Delta_t).
\]
Note that there exists an unique continuous-time Markov process satisfying the above transition rates functions according to classic probability theory (see Section one of \cite{Lig2010}).

Intuitively, $\{X_t\}_{t\geq 0}$ describes the spread of an epidemic on $G_n$. Vertices in state $0$ are susceptible which are healthy and may be infected by the epidemic. Vertices in state $1$ are infective that can infect susceptible neighbors. Vertices in state $-1$ are removed which will never be infected again. An infective vertex waits for an exponential time with rate one to become removed while a susceptible vertex is infected by an infective neighbor at rate proportional to the production of the vertex weights on these two vertices. Note that here we say two vertices are neighbors when the edge connecting them is remained during the procedure to generate $G_n$.

For any $t\geq 0$, we define
\begin{equation}\label{equ 1.1}
S_t=\sum_{j=0}^{n-1}\mathbf{1}_{\{X_t(j)=0\}}
\end{equation}
as the number of susceptible vertices at moment $t$ and
\begin{equation}\label{equ 1.2}
V_t=\sum_{j=0}^{n-1}\rho(j)\mathbf{1}_{\{X_t(j)=1\}}
\end{equation}
as the total capability of the infective vertices to infect neighbors at moment $t$. We write $X_t, ~S_t$ and $V_t$ as $X_t^{(n)}, ~S_t^{(n)}$ and $V_t^{(n)}$ when we need to point out that process is on $G_n$. For the moment $t=0$, we assume that $\{X_0^{(n)}(i)\}_{i=0}^{n-1}$ are i. i. d. such that
\begin{equation}\label{equ 1.3}
P(X_0^{(n)}(0)=1)=\theta=1-P(X_0^{(n)}(0)=0)
\end{equation}
for some $\theta\in (0,1)$. Under this assumption, we obtain the law of large numbers for $\big(\frac{S_t^{(n)}}{n},~\frac{V_t^{(n)}}{n}\big)$ as $n$ grows to infinity at any moment $t$. For mathematical details, see the next section.

Readers may wonder what will occur when at $t=0$ there is only one infective vertex. We study a similar epidemic model under this assumption in \cite{Xue 2016}. According to a similar analysis with that in \cite{Xue2016}, it can be shown that if there is only one infective vertex at $t=0$, then the epidemic `outbreaks' when and only when $\lambda>\lambda_c=\frac{1}{p{\rm E}\rho^2}$ (see the main theorem of \cite{Xue2016} for the accurate meaning of `outbreak'). Actually we think this is an important result but we do not want to repeat in this paper lot of similar calculation with that in \cite{Xue2016}, so we only give a simple comment here. Readers interested with mathematical details can see \cite{Xue2016}.

When $\rho=p=1$, then our model reduces to the classical SIR model on complete graphs. According to the theory of density dependent population model introduced in Section 11 of \cite{Ethier1986}, under assumption \eqref{equ 1.3}, $(\frac{S_t^{(n)}}{n},~\frac{V_t^{(n)}}{n})$ converges in probability to the solution $(s_t,~v_t)$ of the following ODE as $n\rightarrow+\infty$.
\begin{equation}\label{equ 1.4 ODE}
\begin{cases}
&\frac{d}{dt}s_t=-\lambda s_tv_t,\\
&\frac{d}{dt}v_t=\lambda s_tv_t-v_t,\\
&s_0=1-\theta,\\
&v_0=\theta.
\end{cases}
\end{equation}
Our main result given in the next section can be seen as an extension of the above conclusion.

In \cite{Neal2003}, Neal studies a discrete-time version of SIR on Erd\H{o}s-R\'{e}nyi graph $G(n,p)$. The exact final-size distribution and extensive asymptotic results of the model are given since it is shown that the model is a randomized Reed-Frost epidemic, property of which is studied in \cite{Martin1986} and \cite{Picard1990}.

There is another important type of epidemic model which is the SIS model. For the SIS model, an infective one will become susceptible and can be infected again. The stochastic SIS model is also named as the contact process. Readers can see Chapter 6 of \cite{Lig1985} and Part one of \cite{Lig1999} for a survey of the study of the contact process. Peterson studies the contact process with random vertex weights on complete graphs in \cite{Pet2011} and gives the critical value of the process, which inspires us a lot.

Erd\H{o}s-R\'{e}nyi graph can be considered as the bond percolation model on the complete graph (see \cite{Grimmett1999} for the definition of the
percolation model). In recent years, epidemic on clusters of the percolation model has been a popular topic. Readers can see references
\cite{Ber2011, Chen2009, Xue2016b, Xue2016c} and so on for related study.

In \cite{Janson2014}, Janson, Luczak and Windridge study the SIR
epidemic on a random graph with given degrees. We are inspired a lot
by the law of large numbers given in \cite{Janson2014}. Our LLN
given in the next section is with a similar form with that in
\cite{Janson2014}. In \cite{Kurtz2008}, Kurtz and co-workers 
study limit theorems for an epidemic model on the complete graph. We solve an ODE in this paper according to the approach introduced in \cite{Kurtz2008}.

\section{Main result}\label{section two}
In this section we will give our main results. First we introduce
some notations. For any $x>0$, we define
\begin{equation}\label{equ 2.1}
H_S(x)=(1-\theta){\rm E}\big(x^{\rho}\big)
\end{equation}
and
\begin{equation}\label{equ 2.2}
H_V(x)={\rm E}\rho-(1-\theta){\rm E}\big(\rho
x^{\rho}\big)+\frac{\log x}{p\lambda},
\end{equation}
where ${\rm E}$ is the expectation operator with respect to the
random variable $\rho$ while $\theta,p$ and $\lambda$ are introduced
as in Section \ref{section one}. Since $P(0\leq \rho\leq M_1)=1$, it
is easy to check that $H_V(x)$ satisfies the local Lipschitz
condition. As a result, the following ODE
\begin{equation}\label{equ 2.3}
\begin{cases}
&\frac{d}{dt}\psi_t=-p\lambda\psi_tH_V(\psi_t),\\
&\psi_0=1
\end{cases}
\end{equation}
has an unique solution $\{\psi_t\}_{t\geq 0}$. Now we can give our
main result, which is the following law of large numbers.
\begin{theorem}\label{theorem 2.1 main LLN}
If $S_t, ~V_t, ~H_S(x)$ and $H_V(x)$ are defined as in Equations
\eqref{equ 1.1}, \eqref{equ 1.2}, \eqref{equ 2.1} and \eqref{equ
2.2} respectively while $\{\psi_t\}_{t\geq 0}$ is the unique
solution of ODE \eqref{equ 2.3}, then under assumption \eqref{equ
1.3},
\[
\lim_{n\rightarrow+\infty}\frac{S_t^{(n)}}{n}=H_S(\psi_t)\text{\quad
and\quad} \lim_{n\rightarrow+\infty}\frac{V_t^{(n)}}{n}=H_V(\psi_t)
\]
in probability for any $t\geq 0$.
\end{theorem}

Here we give an non-rigorous explanation of Theorem \ref{theorem 2.1
main LLN}. Intuitively, $\psi_t$ is the limit probability that a
given initially susceptible vertex with weight $1$ is still
susceptible at the moment $t$. By a mean-field idea, an susceptible
vertex with weight $\rho$ is infected at rate approximate $\rho$
times that of vertex with weight one. Hence the probability that an
initially susceptible vertex with weight $\rho$ is still susceptible
at the moment $t$ is about $\psi_t^\rho$. Initially susceptible
vertices with weight $\rho$ have density about $(1-\theta)P(d\rho)$
by assumption \eqref{equ 1.3}. Hence, the proportion of $S_t$ is
about
\[
\int \psi_t^\rho(1-\theta)P(d\rho)=(1-\theta){\rm
E}\big(\psi_t^{\rho}\big)=H_S(\psi_t)
\]
for large $n$.

The proof of Theorem \ref{theorem 2.1 main LLN} is given in the next
section, where we will assume that $\rho$ has finite support. This
assumption is without loss of generality according to the following
analysis. For $\rho$ with infinite support and integer $m\geq 1$, we
let $\widehat{\rho}_m=\frac{\lfloor m\rho\rfloor}{m}$ and
$\widetilde{\rho}_m=\frac{\lfloor m\rho\rfloor+1}{m}$, then
$\widehat{\rho}_m$ and $\widetilde{\rho}_{m}$ have finite support
while
\begin{equation}\label{equ 2.4}
\lim_{m\rightarrow+\infty}\widehat{\rho}_m=\lim_{m\rightarrow+\infty}\widetilde{\rho}_m=\rho.
\end{equation}
We define $\widehat{S}_{t,m},\widehat{V}_{t,m}, \widehat{H}_{S,m},
\widehat{H}_{V,m},\widehat{\psi}_{t,m}$ as counterparts with respect
to $\widehat{\rho}_m$ of $S_t,V_t, H_S, H_V, \psi_t$ while
$\widetilde{S}_{t,m},\widetilde{V}_{t,m}, \widetilde{H}_{S,m},
\widetilde{H}_{V,m},\widetilde{\psi}_{t,m}$ as counterparts with
respect to $\widetilde{\rho}_m$ of $S_t,V_t, H_S, H_V,\psi_t$. By
Equation \eqref{equ 2.4},
\begin{align}\label{equ 2.5}
&\lim_{m\rightarrow+\infty}\widetilde{H}_{S,m}(\widetilde{\psi}_{t,m})=\lim_{m\rightarrow+\infty}\widetilde{H}_{S,m}(\widetilde{\psi}_{t,m})
=H_S(\psi_t)\notag\\
\text{~and\quad}
&\lim_{m\rightarrow+\infty}\widetilde{H}_{V,m}(\widetilde{\psi}_{t,m})=\lim_{m\rightarrow+\infty}\widetilde{H}_{V,m}(\widetilde{\psi}_{t,m})
=H_V(\psi_t).
\end{align}
By basic coupling of Markov processes (see Section 2.1 of
\cite{Lig1985}),
\begin{equation}\label{equ 2.6}
\widehat{S}_{t,m}\geq S_t\geq \widetilde{S}_{t,m}\text{~and~}
\widehat{V}_{t,m}\leq V_t\leq \widetilde{V}_{t,m}
\end{equation}
in the sense of coupling, since $\widehat{\rho}_m\leq \rho\leq
\widetilde{\rho}_m$. By Equations \eqref{equ 2.5} and \eqref{equ
2.6}, Theorem \ref{theorem 2.1 main LLN} will hold for $(S_t, V_t)$
if it is proved that Theorem \ref{theorem 2.1 main LLN} holds for
$(\widehat{S}_{t,m},\widehat{V}_{t,m})$ and
$(\widetilde{S}_{t,m},\widetilde{V}_{t,m})$. As a result, we only
need to deal with the case where $\rho$ has finite support.

\section{Proof of Theorem \ref{theorem 2.1 main LLN}}\label{section three}
In this section we will give the proof of Theorem \ref{theorem 2.1 main LLN}. First we introduce some notations and definitions. Without loss of generality, we assume that the vertex weight $\rho$ has finite support as we have introduced at the end of Section \ref{section two}. Let $q_1,q_2,\ldots,q_K\in [0,M_1]$ and $\mu_1,\mu_2,\ldots,\mu_K$ be strictly positive constants such that $\sum\limits_{j=1}^{K}\mu_j=1$, then we assume that
\begin{equation}\label{equ 3.-1}
P(\rho=q_j)=\mu_j
\end{equation}
for $1\leq j\leq K$, where $M_1$ is defined as in Section \ref{section one}. For each $1\leq j\leq K$ and any $t\geq 0$, we define
\[
S_t^{(n)}(j)={\rm card}\big\{0\leq i\leq n-1:~X_t^{(n)}(i)=0 \text{~and~}\rho(i)=q_j\big\}
\]
and
\[
I_t^{(n)}(j)={\rm card}\big\{0\leq i\leq n-1:~X_t^{(n)}(i)=1 \text{~and~}\rho(i)=q_j\big\}
\]
as the numbers of susceptible vertices and infective vertices respectively with weight $q_j$ at the moment $t$, where we use ${\rm card}(\cdot)$ to denote the cardinality of the set. Hence,
\begin{equation}\label{equ 3.0}
S_t^{(n)}=\sum_{j=1}^KS_t^{(n)}(j)\text{~and~} V_t^{(n)}=\sum_{j=1}^K\rho(j)I_t^{(n)}(j).
\end{equation}
We write $S_t^{(n)}(j)$ and $I_t^{(n)}(j)$ as $S_t(j)$ and $I_t(j)$ when there is no misunderstanding.

For $C,D\subseteq A_n$ such that $C\bigcap D=\emptyset$, we define
\[
\alpha(C,D)=\sum_{i\in C,~j\in D}\mathbf{1}_{\{i\sim j\}}
\]
as the number of edges connecting two vertices in $C$ and $D$ respectively. For later use, for any $0<c,d\leq 1$ and each $n\geq 1$, we define
\begin{align}\label{equ 3.1}
\beta(c,d,n)=\sup\Big\{&\big|\alpha(C,D)-{\rm card}(C){\rm card}(D)p\big|:~C\subseteq A_n, D\subseteq A_n, \notag\\
&C\bigcap D=\emptyset, {\rm card}(C)\geq cn, {\rm card}(D)\geq dn\Big\}.
\end{align}

As a preparation of the proof, we give the following two lemmas.
\begin{lemma}\label{lemma 3.1}
For any $t>0$ and $1\leq j\leq K$,
\[
\lim_{n\rightarrow+\infty}P\big(\inf_{0\leq u\leq t}\frac{S_u^{(n)}(j)}{n}~\geq (1-\theta)\mu_je^{-2\lambda M_1^2t}\big)=1
\]
and
\[
\lim_{n\rightarrow+\infty}P\big(\inf_{0\leq u\leq t}\frac{I_u^{(n)}(j)}{n}~\geq \theta\mu_j e^{-2t}\big)=1.
\]
\end{lemma}

\begin{lemma}\label{lemma 3.2}
For any $c,d>0$ such that $c+d\leq 1$,
\[
\lim_{n\rightarrow+\infty}\frac{\beta(c,d,n)}{n^2}=0
\]
in probability.
\end{lemma}

\proof[Proof of Lemma \ref{lemma 3.1}]

According to the definition of our model on $G_n$, a susceptible vertex is infected at rate at most
\[
n\frac{\lambda}{n}M_1^2=\lambda M_1^2,
\]
according to which we consider the following auxiliary model. We assume that at $t=0$, there are $n$ spins taking value from $\{2,-2\}$. Each spin independently takes $2$ with probability $(1-\theta)\mu_j$ while takes $-2$ with probability $1-(1-\theta)\mu_j$. A spin in state $-2$ will be frozen in this state forever while a spin in state $2$ waits for an exponential time with rate $\lambda M_1^2$ to become $-2$ and then be frozen in state $-2$. Let $Y_t^{(n)}$ be the number of spins in state $2$ at the moment $t$, then $S_0^{(n)}(j)=Y_0^{(n)}$ and $S_t^{(n)}(j)\geq Y_t^{(n)}$ for any $t\geq 0$ in the sense of coupling. By Kolmogorov's law of large numbers,
\begin{equation}\label{equ 3.2}
\lim_{n\rightarrow+\infty}\sup_{0\leq u\leq t}\big|\frac{Y_u^{(n)}}{n}-(1-\theta)\mu_je^{-\lambda M_1^2u}\big|~=0
\end{equation}
in probability since a given spin is in state $2$ at moment $u$ with probability $(1-\theta)\mu_je^{-\lambda M_1^2u}$. The first part of Lemma 3.1 follows directly from Equation \eqref{equ 3.2} and the fact that $S_t^{(n)}(j)\geq Y_t^{(n)}$ for any $t\geq0$.

The proof of the second part is similar. We consider another $n$ spins taking value from $\{3,-3\}$ and assume that a spin in state $-3$ will be frozen in this state forever while a spin in state $3$ becomes $-3$ at rate one and then is frozen in state $-3$. We denote by $Z_t^{(n)}$ the number of spins in state $3$ at moment $t$ and assume that $Z_0^{(n)}=I_0^{(n)}(j)$, then $I_t^{(n)}(j)\geq Z_t^{(n)}$ for any $t\geq 0$ in the sense of coupling since any infective vertex becomes removed at rate one but infective vertices may increase through infecting susceptible ones. The second part of Lemma \ref{lemma 3.1} follows from this coupling and the fact that
\[
\lim_{n\rightarrow+\infty}\sup_{0\leq u\leq t}\big|\frac{Z_u^{(n)}}{n}-\theta\mu_je^{-u}\big|~=0
\]
in probability according to Kolmogorov's law of large numbers.

\qed

\proof[Proof of Lemma \ref{lemma 3.2}]

Let $\{U_j\}_{j=1}^{+\infty}$ be i. i. d. random variables such that
\[
P(U_1=1)=p=1-P(U_1=0),
\]
then according to classic theory of large deviation principle, for any $\epsilon>0$ and integer $N\geq 1$,
\begin{equation}\label{equ 3.3}
P\Big(\big|\frac{\sum_{j=1}^NU_j}{N}-p\big|\geq \epsilon\Big)\leq \big[\sigma(\epsilon)\big]^N,
\end{equation}
where $\sigma(\epsilon)\in (0,1)$. By Equation \eqref{equ 3.3}, for $C,D\subseteq A_n$ such that $C\bigcap D=\emptyset$, ${\rm card}(C)\geq cn$ and ${\rm card}(D)\geq dn$,
\begin{align}\label{equ 3.4}
&P\Big(\big|\alpha(C,D)-{\rm card}(C){\rm card}(D)p\big|\geq n^2\epsilon\Big)\notag\\
&\leq P\Big(\big|\frac{\alpha(C,D)}{{\rm card}(C){\rm card}(D)}-p\big|\geq \epsilon\Big)\leq \big[\sigma(\epsilon)\big]^{{\rm card}(C){\rm card}(D)}\\
&\leq \big[\sigma(\epsilon)\big]^{cdn^2},\notag
\end{align}
since $\alpha(C,D)$ has the same probability distribution as that of $\sum\limits_{j\leq {\rm card}(C){\rm card}(D)}U_j$. The number of subsets of a set with $n$ elements is $2^n$, hence by Equation \eqref{equ 3.4},
\begin{equation}\label{equ 3.5}
P(\beta(c,d,n)\geq n^2\epsilon)\leq 2^n2^n \big[\sigma(\epsilon)\big]^{cdn^2}=4^n \big[\sigma(\epsilon)\big]^{cdn^2}\rightarrow 0
\end{equation}
as $n\rightarrow+\infty$. Lemma \ref{lemma 3.2} follows from Equation \eqref{equ 3.5} directly. 

\qed

For $1\leq j,l\leq K$ and $t\geq 0$, we denote by $L_t^{(n)}(j,l)$ the number of edges connecting a susceptible vertex with weight $q_j$ and an infective vertex with weight $q_l$ at moment $t$. The superscript $(n)$ means that the process is on $G_n$. According to Lemmas \ref{lemma 3.1} and \ref{lemma 3.2}, we have the following important corollary about $L_t^{(n)}$.

\begin{corollary}\label{corollary 3.3}
For any $1\leq j,l\leq K$ and $t\geq 0$,
\[
\lim_{n\rightarrow+\infty}\frac{\sup\limits_{0\leq u\leq t}\big|L_u^{(n)}(j,l)-pS_u^{(n)}(j)~I_u^{(n)}(l)\big|}{n^2}~=0
\]
in probability.
\end{corollary}

\proof

By Lemma \ref{lemma 3.1}, except for a set with arbitrarily small probability as $n\rightarrow+\infty$,
\begin{equation}\label{equ 3.6}
\sup\limits_{0\leq u\leq t}\big|L_u^{(n)}(j,l)-pS_u^{(n)}(j)~I_u^{(n)}(l)\big|\leq \beta\Big( (1-\theta)\mu_je^{-2\lambda M_1^2t},~\theta\mu_j e^{-2t},~n\Big).
\end{equation}
Corollary \ref{corollary 3.3} follows from Equation \eqref{equ 3.6} and Lemma \ref{lemma 3.2} directly.

\qed

According to corollary \ref{corollary 3.3}, we obtain the following law of large numbers of our model through a similar analysis with that Ethier and Kurtz utilize to construct the theory of density dependent population model in \cite{Ethier1986}.

\begin{proposition}\label{proposition 3.4}
For any $t\geq 0$,
$\Big(\frac{S_t^{(n)}(1)}{n},\frac{S_t^{(n)}(2)}{n},\ldots, \frac{S_t^{(n)}(K)}{n},\frac{V_t^{(n)}}{n}\Big)$ converges in probability to the solution $\big(s_t(1),s_t(2),\ldots,s_t(K),v_t\big)$ of the ODE
\begin{equation}\label{equ 3.7 ODE}
\begin{cases}
&\frac{d}{dt}s_t(i)=-p\lambda v_t q_i s_t(i) \text{\quad for~}1\leq i\leq K,\\
&\frac{d}{dt}v_t=-v_t+p\lambda v_t\big(\sum_{i=1}^Kq_i^2s_t(i)\big),\\
&s_0(i)=(1-\theta)\mu_i \text{\quad for~}1\leq i\leq K,\\
&v_0=\theta \sum_{i=1}^Kq_i\mu_i=\theta {\rm E}\rho
\end{cases}
\end{equation}
as $n\rightarrow+\infty$.
\end{proposition}
Note that it is easy to check that ODE \eqref{equ 3.7 ODE} satisfies local Lipschitz condition, hence the solution $\big(s_t(1),s_t(2),\ldots,s_t(K),v_t\big)$ is unique.

\proof[Proof of Proposition \ref{proposition 3.4}]

According to the evolution of our model, for each $1\leq i\leq K$, $S_t^{(n)}(i)$ decreases by one at rate
\[
\frac{\lambda}{n}q_i\sum_{j=1}^Kq_jL_t^{(n)}(i,j).
\]
Then, by Proposition 1.7 of \cite{Ethier1986}, we can write $S_t^{(n)}(i)$ as
\begin{equation}\label{equ 3.8}
S_t^{(n)}(i)=S_0^{(n)}(i)-\int_0^t\frac{\lambda q_i}{n}\sum_{j=1}^Kq_jL_u^{(n)}(i,j)~du-\widetilde{N}_i\Big(\int_0^t\frac{\lambda q_i}{n}\sum_{j=1}^Kq_jL_u^{(n)}(i,j)~du\Big),
\end{equation}
where $\{\widetilde{N}_i(s):~s\geq 0\}$ is a right-continuous martingale such that ${\rm E}\widetilde{N}_i(s)=0$ and ${\rm E}\big(\widetilde{N}_i^2(s)\big)=s$ for any $s>0$. Furthermore, $\{\widetilde{N}_i\}_{i=1}^K$ are i. i. d..

By Equations \eqref{equ 3.0} and \eqref{equ 3.8}, we can write $S_t^{(n)}(i)$ as
\begin{equation}\label{equ 3.10}
S_t^{(n)}(i)=S_0^{(n)}(i)-\int_0^t\frac{p\lambda q_i S_u^{(n)}(i)V_u^{(n)}}{n}~du+n\epsilon_{1,i}^{(n)}(t)+n\epsilon_{2,i}^{(n)}(t),
\end{equation}
where
\[
\epsilon_{1,i}^{(n)}(t)=-\frac{1}{n}\widetilde{N}_i\Big(\int_0^t\frac{\lambda q_i}{n}\sum_{j=1}^Kq_jL_u^{(n)}(i,j)~du\Big)
\]
and
\[
\epsilon_{2,i}^{(n)}(t)=-\int_0^t\sum_{j=1}^K\lambda q_iq_j\frac{\Big(L_u^{(n)}(i,j)-pS_u^{(n)}(i)I_u^{(n)}(j)\Big)}{n^2}~du.
\]
Since $\big(s_t(1),\ldots,s_t(K),v_t\big)$ is the solution of ODE \eqref{equ 3.7 ODE}, by \eqref{equ 3.10},
\begin{align}\label{equ 3.11}
\frac{S_t^{(n)}(i)}{n}-s_t(i)=&\Big(\frac{S_0^{(n)}(i)}{n}-(1-\theta)\mu_i\Big)-\int_0^tp\lambda q_i \Big(\frac{S_u^{(n)}(i)}{n}\frac{V_u^{(n)}}{n}-s_u(i)v_u\Big)~du  \\
&+\epsilon_{1,i}^{(n)}(t)+\epsilon_{2,i}^{(n)}(t).  \notag
\end{align}
According to the evolution of our model, for each $1\leq i\leq K$, $V_t^{(n)}$ decreases by $q_i$ at rate $I^{(n)}_t(i)$ and increases by $q_i$ at rate
\[
\frac{\lambda}{n}q_i\sum_{j=1}^Kq_jL_t^{(n)}(i,j).
\]
Therefore, by Proposition 1.7 of \cite{Ethier1986} and an similar analysis with that deduces Equation \eqref{equ 3.11},
\begin{align}\label{equ 3.12}
\frac{V_t^{(n)}}{n}-v_t=&\Big(\frac{V_0^{(n)}}{n}-\theta{\rm E}\rho\Big)-\int_0^t\Big(\frac{V_u^{(n)}}{n}-v_u\Big)~du\\
&+\int_{0}^tp\lambda\sum_{i=1}^{K}q_i^2\Big(\frac{S_u^{(n)}(i)}{n}\frac{V_u^{(n)}}{n}-s_u(i)v_u\Big)~du +\epsilon^{(n)}_3(t)+\epsilon^{(n)}_4(t)+\epsilon^{(n)}_5(t),\notag
\end{align}
where
\begin{align*}
\epsilon^{(n)}_3(t)&=\sum_{i=1}^Kq_i^2\int_0^t\frac{\lambda\sum_{j=1}^K q_j\Big(L_u^{(n)}(i,j)-pS_u^{(n)}(i)I_u^{(n)}(j)\Big)}{n^2}~du,\\
\epsilon_4^{(n)}(t)&=-\frac{1}{n}\sum_{i=1}^Kq_i\widetilde{N}_{i,-}\Big(\int_0^tI^{(n)}_u(i)~du\Big)\\
\text{~and\quad}\epsilon_5^{(n)}(t)&=\frac{1}{n}\sum_{i=1}^Kq_i\widetilde{N}_i\Big(\int_0^t\frac{\lambda}{n}q_i\sum_{j=1}^Kq_jL_u^{(n)}(i,j)~du\Big),
\end{align*}
where $\{\widetilde{N}_{i,-}\}_{i=1}^K$ is an independent copy of $\{\widetilde{N}_i\}_{i=1}^K$. Note that since $S_t^{(n)}(i)$ decreasing by one and $V_t^{(n)}$ increasing by $q_i$ occur at the same moment, the martingale with respect to the jumps of $V_t^{(n)}$ increasing by $q_i$ is the same as that with respect to the jumps of $S_t^{(n)}$ decreasing by one, which is $\widetilde{N}_i$. This is why we do not need another independent copy of $\{\widetilde{N}_i\}_{i=1}^K$ in the form of $\epsilon_5$.

For given $T>0$ and $1\leq i\leq K$, we define
\[
\epsilon_6^{(n)}(i)=\frac{1}{n}\sup\Big\{\big|\widetilde{N}_i(s)\big|:~0\leq s\leq \lambda TM_1^2Kn\Big\}
\]
and
\[
\epsilon_7^{(n)}(i)=\frac{1}{n}\sup\Big\{\big|\widetilde{N}_{i,-}(s)\big|:~0\leq s\leq nT\Big\},
\]
then, for any $t\leq T$,
\begin{align}\label{equ 3.13}
&|\epsilon_{1,i}^{(n)}(t)|\leq \epsilon_6^{(n)}(i),~|\epsilon_5^{(n)}(t)|\leq \sum_{i=1}^Kq_i\epsilon_6^{(n)}(i)\\
\text{~and\quad} &|\epsilon_4^{(n)}(t)|\leq \sum_{i=1}^Kq_i\epsilon_7^{(n)}(i) \notag
\end{align}
since
\[
\int_0^t\frac{\lambda q_i}{n}\sum_{j=1}^Kq_jL_u^{(n)}(i,j)~du\leq T\frac{\lambda}{n}KM_1^2n^2=\lambda TM_1^2Kn
\]
and \[
\int_0^tI^{(n)}_u(i)~du\leq nT.
\]
By Doob's inequality,
\begin{align*}
P(\epsilon_6^{(n)}(i)\geq \epsilon)\leq \frac{4}{n^2\epsilon^2}{\rm E}\Big(\widetilde{N}^2_i\big(\lambda TM_1^2Kn\big)\Big)=\frac{4\lambda TM_1^2K}{n\epsilon^2}\rightarrow 0
\end{align*}
as $n\rightarrow+\infty$. Therefore
\begin{equation}\label{equ 3.14}
\lim_{n\rightarrow+\infty}\epsilon_6^{(n)}(i)=0
\end{equation}
in probability. According to a similar analysis,
\begin{equation}\label{equ 3.15}
\lim_{n\rightarrow+\infty}\epsilon_7^{(n)}(i)=0
\end{equation}
in probability. For given $T>0$, we define
\[
\epsilon_8^{(n)}=\frac{\sup\limits_{1\leq i,j\leq K}\sup\limits_{0\leq u\leq T} \Big|L_u^{(n)}(i,j)-pS_u^{(n)}(i)I_u^{(n)}(j)\Big|}{n^2},
\]
then
\begin{equation}\label{equ 3.16}
\lim_{n\rightarrow+\infty}\epsilon_8^{(n)}=0
\end{equation}
in probability according to Corollary \ref{corollary 3.3} and
\begin{equation}\label{equ 3.17}
|\epsilon_{2,i}(t)|\leq \lambda M_1^2TK\epsilon_8^{(n)}
\end{equation}
for $t\leq T$ and $1\leq i\leq K$ while
\begin{equation}\label{equ 3.18}
|\epsilon_3^{(n)}(t)|\leq \lambda K^2M_1^3T\epsilon_8^{(n)}
\end{equation}
for $t\leq T$ since $P(\rho\leq M_1)=1$. For each $t\geq 0$, we define
\[
\sigma_t^{(n)}=\big|\frac{V_t^{(n)}}{n}-v_t\big|+\sum_{i=1}^K\big|\frac{S_t^{(n)}(i)}{n}-s_t(i)\big|.
\]
For given $T>0$ and $1\leq i\leq K$, we define
\[
D_i=\max\big\{1, ~\sup_{0\leq t\leq T}|s_t(i)|\big\}
\]
and $D_{K+1}=\max\big\{M_1,~\sup_{0\leq t\leq T}|v_t|\big\}$, then it is easy to check that there exists $M_2>0$ such that
\begin{equation}\label{equ 3.19}
\sum_{i=1}^K\big|x_iv_1-y_iv_2\big|\leq M_2\Big(\big|v_1-v_2\big|+\sum_{i=1}^K\big|x_i-y_i\big|\Big)
\end{equation}
for any $\vec{x}=(x_1,\ldots,x_K,v_1)\in \prod\limits_{i=1}^{K+1}[-D_i,D_i]$ and $\vec{y}=(y_1,\ldots,y_K,v_2)\in \prod\limits_{i=1}^{K+1}[-D_i,D_i]$.
By Equation \eqref{equ 3.19},
\begin{equation}\label{equ 3.20}
\sum_{i=1}^K\big|\frac{S_u^{(n)}(i)}{n}\frac{V_u^{(n)}}{n}-s_u(i)v_u\big|\leq M_2 \sigma_u^{(n)}
\end{equation}
for any $u\leq T$, since $\frac{S_u^{(n)}}{n}\leq 1$ and $\frac{V_u^{(n)}}{n}\leq M_1$. By Equations \eqref{equ 3.12}, \eqref{equ 3.13}, \eqref{equ 3.18} and \eqref{equ 3.20}, for $t\leq T$,
\begin{equation}\label{equ 3.21}
\big|\frac{V_t^{(n)}}{n}-v_t\big|\leq \epsilon_9^{(n)} +\int_0^t\big|\frac{V_u^{(n)}}{n}-v_u\big|~du+p\lambda M_1^2M_2\int_0^t\sigma_u^{(n)}~du,
\end{equation}
where
\[
\epsilon_9^{(n)}=\big|\frac{V_0^{(n)}}{n}-\theta{\rm E}\rho\big|+\lambda K^2M_1^3T\epsilon_8^{(n)}+\sum_{i=1}^Kq_i\epsilon_6^{(n)}(i)
+\sum_{i=1}^Kq_i\epsilon_7^{(n)}(i).
\]
By Equations \eqref{equ 3.11}, \eqref{equ 3.13}, \eqref{equ 3.17} and \eqref{equ 3.20}, for $t\leq T$,
\begin{equation}\label{equ 3.22}
\sum_{i=1}^K\big|\frac{S_t^{(n)}(i)}{n}-s_t(i)\big|\leq \epsilon_{10}^{(n)}+p\lambda M_1M_2\int_0^t\sigma^{(n)}_u~du,
\end{equation}
where
\[
\epsilon_{10}^{(n)}=\sum_{i=1}^K\big|\frac{S_0^{(n)}(i)}{n}-(1-\theta)\mu_i\big|+\sum_{i=1}^K\epsilon_6^{(n)}(i)+\lambda M_1^2TK^2\epsilon_8^{(n)}.
\]
By Equations \eqref{equ 3.21} and \eqref{equ 3.22},
\begin{equation}\label{equ 3.23}
\sigma_t^{(n)}\leq \epsilon_9^{(n)}+\epsilon_{10}^{(n)}+\big(1+p\lambda M_1^2M_2+p\lambda M_1M_2\big)\int_0^t\sigma_u^{(n)}~du
\end{equation}
for $0\leq t\leq T$. Let $M_3=1+p\lambda M_1^2M_2+p\lambda M_1M_2$, then by Equation \eqref{equ 3.23} and Gronwall's inequality,
\begin{equation}\label{equ 3.24}
\sigma_t^{(n)}\leq \big(\epsilon_9^{(n)}+\epsilon_{10}^{(n)}\big)e^{M_3t}.
\end{equation}
for $0\leq t\leq T$. By assumptions \eqref{equ 1.3}, \eqref{equ 3.-1} and Kolmogorov's law of large numbers,
\[
\lim_{n\rightarrow+\infty}\big|\frac{V_0^{(n)}}{n}-\theta{\rm E}\rho\big|+\sum_{i=1}^K\big|\frac{S_0^{(n)}(i)}{n}-(1-\theta)\mu_i\big|=0
\]
in probability. Therefore,
\begin{equation}\label{equ 3.25}
\lim_{n\rightarrow+\infty}\big(\epsilon_9^{(n)}+\epsilon_{10}^{(n)}\big)=0
\end{equation}
in probability according to Equations \eqref{equ 3.14}, \eqref{equ 3.15} and \eqref{equ 3.16}. By Equations \eqref{equ 3.24} and \eqref{equ 3.25},
\[
\lim_{n\rightarrow+\infty}\sigma_t^{(n)}=0
\]
in probability for any $0\leq t\leq T$. Since we can choose $T$ arbitrarily large,
\begin{equation}\label{equ 3.27}
\lim_{n\rightarrow+\infty}\sigma_t^{(n)}=0
\end{equation}
in probability for any $t\geq 0$. Note that $\sigma_t^{(n)}=\big|\frac{V_t^{(n)}}{n}-v_t\big|+\sum_{i=1}^K\big|\frac{S_t^{(n)}(i)}{n}-s_t(i)\big|$, hence Proposition \ref{proposition 3.4} follows from Equation \eqref{equ 3.27} directly.

\qed

At last we give the proof of Theorem \ref{theorem 2.1 main LLN}.

\proof[Proof of Theorem \ref{theorem 2.1 main LLN}] In this proof we write $v_t$ as $v(t)$ and write $s_t(i)$ as $s(t,i)$ to make them look better.
According to Proposition \ref{proposition 3.4} and Equation \eqref{equ 3.0}, we only need to show that
\begin{equation}\label{equ 3.28}
v(t)=H_V(\psi_t)\text{~and~} \sum_{i=1}^{K}s(t,i)=H_S(\psi_t),
\end{equation}
where $\psi_t$ is the solution of ODE \eqref{equ 2.3}. We solve ODE \eqref{equ 3.7 ODE} through an approach which is introduced by Kurtz and co-workers in \cite{Kurtz2008}. We define $\{A_t\}_{t\geq 0}$ as the solution of the ODE
\begin{equation}\label{equ 3.29}
\begin{cases}
&\frac{d}{dt}A_t=\frac{1}{v(A_t)},\\
&A_0=0.
\end{cases}
\end{equation}
Then by the definition of $\big(s(t,1),\ldots,s(t,K),v(t)\big)$,
\[
\begin{cases}
&\frac{d}{dt}s(A_t,i)=-p\lambda q_i s(A_t,i) \text{\quad for~}1\leq i\leq K,\\
&\frac{d}{dt}v(A_t)=-1+p\lambda \big(\sum_{i=1}^Kq_i^2s(A_t,i)\big),\\
&s(A_0,i)=(1-\theta)\mu_i \text{\quad for~}1\leq i\leq K,\\
&v(A_0)=\theta \sum_{i=1}^Kq_i\mu_i=\theta {\rm E}\rho.
\end{cases}
\]
Then it is easy to see that
\[
\begin{cases}
s(A_t,i)=(1-\theta)\mu_ie^{-p\lambda q_it}\text{\quad for~}1\leq i\leq K,\\
v(A_t)={\rm E}\rho-t-(1-\theta)\sum_{i=1}\mu_iq_ie^{-p\lambda q_i t}
\end{cases}
\]
and hence
\begin{equation}\label{equ 3.30}
\begin{cases}
&v(A_t)=H_V(e^{-p\lambda t}),\\
&\sum_{i=1}^{K}s(A_t,i)=H_S(e^{-p\lambda t})
\end{cases}
\end{equation}
according to the definitions of $H_S$ and $H_V$ given by Equations \eqref{equ 2.1} and \eqref{equ 2.2}. Let $\phi_t$ be the inverse function of $A_t$ such that
$A_{\phi_t}=t$, then by Equation \eqref{equ 3.30},
\begin{equation}\label{equ 3.31}
v(t)=H_V(e^{-p\lambda \phi_t}) \text{~and~} \sum_{i=1}^{K}s(t,i)=H_S(e^{-p\lambda \phi_t}).
\end{equation}
At last we only need to show that $\psi_t:=e^{-p\lambda \phi_t}$ is the solution of ODE \eqref{equ 2.3}. By the definition of $\phi_t$ and $A_t$,
\begin{align*}
\frac{d}{dt}\psi_t&=-p\lambda e^{-p\lambda \phi_t}\frac{d}{dt}\phi_t=-p\lambda \psi_t\frac{1}{\frac{dA_s}{ds}|_{s=\phi_t}}\\
&=-p\lambda \psi_tv(A_s)|_{s=\phi_t}=-p\lambda \psi_tH_V(e^{-p\lambda \phi_t})=-p\lambda\psi_t H_V(\psi_t).
\end{align*} 
Note that $\phi_0=0$ and $e^{-p\lambda \phi_0}=1$ since $A_0=0$, hence $\psi_t=e^{-p\lambda \phi_t}$ is the solution of ODE \eqref{equ 2.3} and the proof is complete.  

\qed

\quad

\textbf{Acknowledgments.} The author is grateful to the financial
support from the National Natural Science Foundation of China with
grant number 11501542.

{}
\end{document}